\documentclass[11pt,a4paper]{article}
\usepackage{graphicx} 
\usepackage{url} 

\title{A Quasiperiodic Tiling With 12-Fold Rotational Symmetry and Inflation 
	Factor $1+\sqrt{3}$}

\author{Theo P. Schaad\textsuperscript{1} and Peter Stampfli\textsuperscript{2}
	\vspace{10pt}\\
\textsuperscript{1} 2704 38th Ave SW, Seattle, WA 98126, USA; theoschaad@gmail.com\\
\textsuperscript{2} Rue de Lausanne 1, 1580 Avenches, Switzerland; pestampf@gmail.com\\
}

\date{}	
\begin{document}

\maketitle

\begin{abstract}
	
We show how we found substitution rules for a quasiperiodic tiling with local 12-fold rotational symmetry and inflation factor $1 + \sqrt{3}$. The base tiles are a square, a rhomb with an acute angle of 30 degrees, and equilateral triangles that are cut in half. These half-triangles follow three different substitution rules and can be recombined into equilateral triangles in nine different ways to make minor variations of the tiling. The tiling contains quasiperiodically repeated 12-fold rosettes. A central rosette can be enlarged to make an arbitrarily large tiling with 12-fold rotational symmetry. An online computer program is provided that allows the user to explore the tiling\cite{app}. 
	
\end{abstract}

\section*{About Tilings With 12-Fold Rotational Symmetry}
We use the substitution method for creating arbitrarily large tilings. It consists of finding base tiles that can be enlarged and then subdivided with the same base tiles without any gaps or overlaps. The enlargement of their edges is called the inflation factor. Each repetition of enlargement and subdivision creates a new generation of tiles. For finite patches of a tiling this increases the number of tiles, whereas an infinite tiling remains unchanged.

\begin{figure}
	\centering
	\begin{minipage}[b]{1\textwidth}
		\textit{
			\centering
			\includegraphics[width=0.9\linewidth]{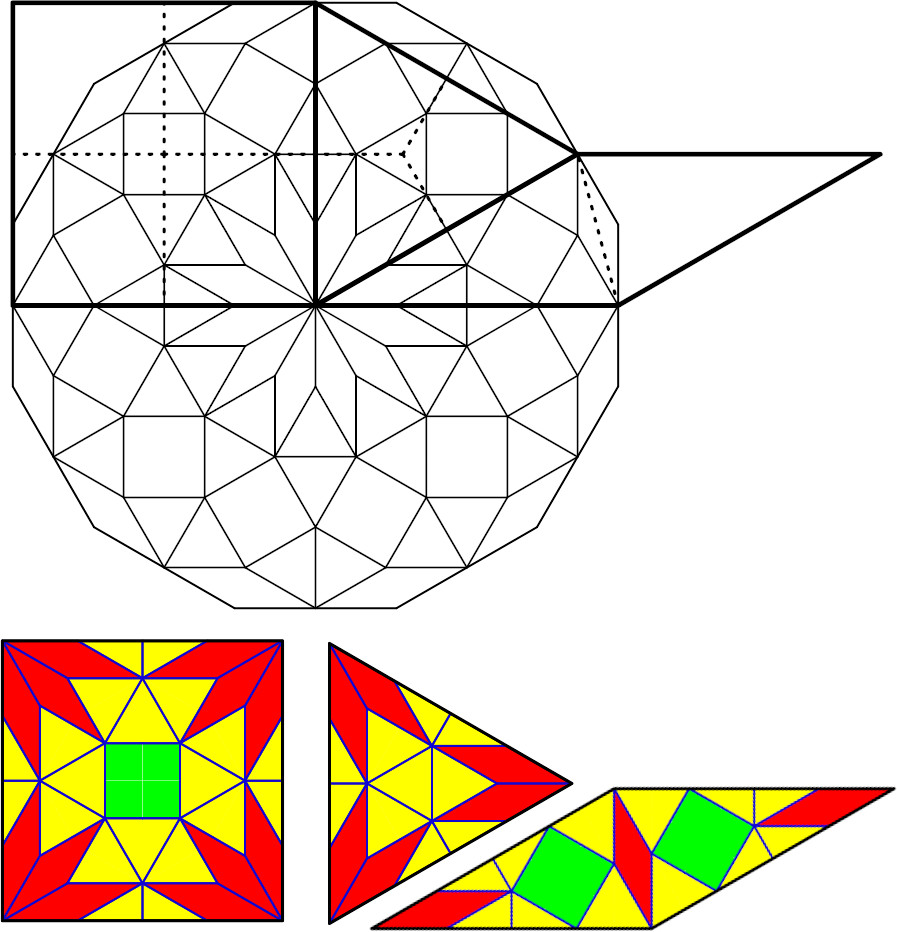}
			\caption{A dodecagon made of two different rhombs and squares. The 12 skinny rhombs form a rosette with local 12-fold rotational symmetry. The rhombs with $60^\circ$ acute angles are further divided into two equilateral triangles, foreshadowing the role it has in this paper. Some quite natural subdivisions of tiles in the dodecagon are shown. They define a substitution scheme for squares, rhombs, and triangles that are inflated by $2 + \sqrt{3}$. All sides of the inflated tiles have the same mirror-symmetric geometry and can be attached to each other.
			}
			\label{fig:rosette}
		}
	\end{minipage}
\end{figure}

Since 1973, many quasiperiodic tilings with $n$-fold rotational symmetry have been discovered, of which the 5-fold Penrose and the 8-fold Ammann-Beenker tilings were the first\cite{gruenbaum}. We consider a special 12-fold quasiperiodic tiling. It exhibits local 12-fold rotational symmetry that is repeated throughout the tiling without any periodicity. Peter Stampfli came up with a particularly simple substitution rule that he described on his blog "Geometry in Color"\cite{stampfli2} and is also listed in the Tiling Encyclopedia of the University of Bielefeld\cite{encyclopedia}. It serves as an introduction here and shows that it is much easier to find a substitution scheme with an inflation factor of $2 + \sqrt{3}$ than for $1 + \sqrt{3}$, which is the subject of this paper. A good starting point for any $n$-fold tiling is an $n$-fold star or rosette of rhombs with acute angles of $2\pi/n$ as shown in the upper part of Figure~\ref{fig:rosette}. The star can then be surrounded with other tiles until an $n$-fold regular polygon is obtained, in this case a dodecagon. In our example, this requires two different rhombs and a square. Peter Stampfli further subdivided the rhombs with $60^{\circ}$ angles into two equilateral triangles. He considered the resulting three shapes to be base tiles of the tiling, namely, a skinny rhomb, a square, and an equilateral triangle. Some quite natural subdivisions within the dodecagon of Figure~\ref{fig:rosette} lead to a substitution scheme with an inflation factor of $2 + \sqrt{3}$. The edges comprise two side lengths of the rhombs and two heights of equilateral triangles. They are conveniently laid out such that all sides of the inflated tiles have the same mirror symmetric pattern of halves of equilateral triangles. These halves meet whenever two tiles are joined edge-to-edge. Thus, any layout using these base tiles can be inflated into an arbitrarily large nonperiodic tiling of equilateral triangles, rhombs and squares. 

Using halves of equilateral triangles has many advantages. This is more efficient and each iteration gives a patch of the same shape as before, except for inflation. It is nice to be able to create square and triangular patches. Using instead full equilateral triangles would give patches that have fractal borders similar to Koch snowflakes because both result from similar iterative methods. Iteration schemes with larger inflation ratios $n+m*\sqrt{3}$, $n\ge 2$, are just as easy to find. If $n$ and $m$ are both odd then we have to split squares at the tile edges.

\begin{figure}
	\centering
	\begin{minipage}[b]{0.9\textwidth}
		\textit{
			\centering
			\includegraphics[width=0.4\linewidth]{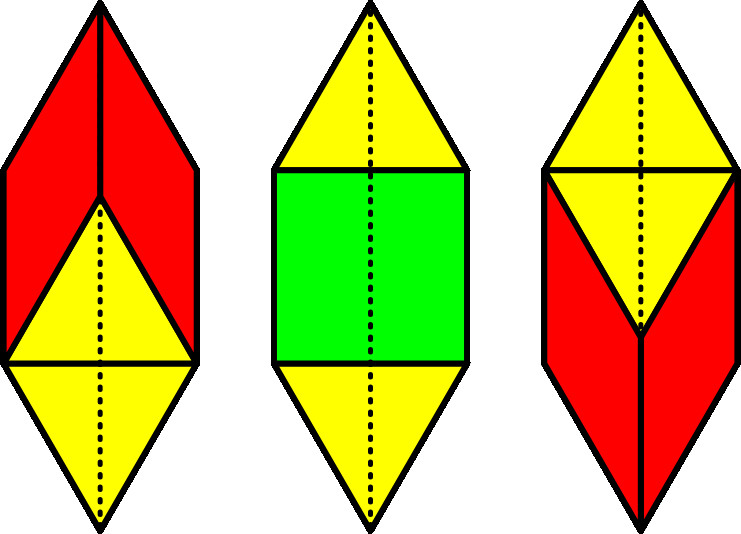}
			\caption{Different dissections of an irregular hexagon into triangles and rhombs or triangles and squares. They can be exchanged in substitution rules, resulting in different tilings. Further, these are all possible substitutions at the sides of tiles (dotted lines) for a self-similarity ratio of $1+\sqrt{3}$.}
			\label{fig:replacement}
		}
	\end{minipage}
\end{figure}

Peter Stampfli further noted that the combination of a triangle with a square can be substituted by another triangle with two skinny rhombs as shown in Figure~\ref{fig:replacement}. This leads to many tilings that sometimes look drastically different\cite{stampfli2}, even though they are only variations. In some cases, the rhombs or the squares can be eliminated and the appearance of local 12-fold rotational symmetry as well\cite{geometricolor}. The Tiling Encyclopedia\cite{encyclopedia} contains several tilings with 12-fold rotational symmetry and inflation factor $2 + \sqrt{3}$. Such tilings are presented by J.~Socolar\cite{socolarOctagonal}, M.~Schlottmann\cite{periodicPointSet}, and Y.~Watanabe et al\cite{watanabe}. J.~Socolar’s first tiling was discovered in 1987 and used regular hexagons, which can also be viewed as three rhombs. This pattern exhibits dodecagonal features, such as dodecagons filled with other tiles, but no outright 12-fold rotational symmetry. The “Shield” tiling by F.~Gähler\cite{gaehlerCrstallography} has the smallest inflation factor of $(2 + \sqrt{3})^{1/2}$. The base tiles consist of an irregular hexagon, a square, and two equilateral triangles with different matching rules. A larger patch exhibits dodecagonal features but no local 12-fold symmetries.

The ”Rorschach” tiling discovered by D.~Frettlöh\cite{rorschach,singThesis} has great relevance to our 12-fold tiling. The inflation factor $1 + \sqrt{3}$ is the same and it uses similar base tiles: a rhomb with 12-fold dihedral symmetry (with an acute angle of $2 \pi/12$), a square, and two equilateral triangles with different substitution rules. The substitution rules prohibit the formation of a 12-fold rosette with rhombs and the resulting tiling has no 12-fold local symmetries. The result is quite surprising and is decidedly a quasiperiodic 4-fold tiling, unique and original in its own rights.

There are other methods than substitution to create dodecagonal tilings. They are only briefly mentioned here, and the reader is referred to the literature. R. Ammann\cite{gruenbaum} discovered that rhombic tilings can be decorated with straight lines that form a grid. This grid is by itself a quasiperiodic tiling called the dual. The dual of this grid is again the original tiling. For the dodecagonal tiling the grid has six sets of parallel lines of different orientation, also called a hexa-grid of parallel lines. It was first explored by J.~Socolar\cite{sololarThesis} in 1987.

A year earlier, Peter Stampfli found another method for creating a dodecagonal tiling\cite{stampfliHelvPhys}. His grid uses periodic tilings of hexagons like a honeycomb. They have a periodic tiling of equilateral triangles as their dual. He discovered that two hexagon grids rotated by 30 degrees to each other have as dual a 12-fold tiling with skinny rhombs, squares, and equilateral triangles. However, this transformation required fairly complicated mathematics. By 2012, he not only posted the results on his blog\cite{geometricolor}, but he also made the computer code freely available. It took several more years to disseminate, but it made the front cover of Science\cite{altourian} in 2018 to illustrate Graphene Quasicrystals made by twisted bilayers of the hexagonal graphene. This 12-fold quasiperiodic tiling does not lend itself to a substitution scheme, but it spurred his interest into finding others.

The work of J.~Socolar was generalized in 1994 by D.~Haussler, H.U.~Nissen, and R.~Lueck\cite{hausslerNissen} in a paper called "Dodecagonal Tilings Derived as Duals from Quasiperiodic Ammann-Grids". They discovered 12 tilings made of rhombs and regular or irregular hexagons but without triangles. The rhombs can be decorated to make new magnified copies of themselves, all with a characteristic deflation factor of $2 + \sqrt{3}$. The results cannot be reproduced by substitution alone. W.~Steurer and S.~Deloudi\cite{stuererDeludy} also discuss the grid method and show that a tiling with a deflation ratio of $1 + \sqrt{3}$ can be obtained. However, it does not have distinct triangles as tiles and there are only few rosettes of rhombs.

To find new tilings we search for new substitution rules. This is an interesting and challenging geometric puzzle similar to tangram and dissection puzzles. We use different approaches for trying out substitution rules and to see if they result in a large tiling without gaps, overlaps, or broken tiles. Theo Schaad draws the substitution rules on paper, scans and copies them, then cuts them out and pastes them together to get higher iterations. Peter Stampfli programs the rules on a computer and lets it iterate them. Both can see if the substitution rules work. 

\begin{figure}
	\centering
	\begin{minipage}[b]{1\textwidth}
		\textit{
			\centering
			\includegraphics[width=0.9\linewidth]{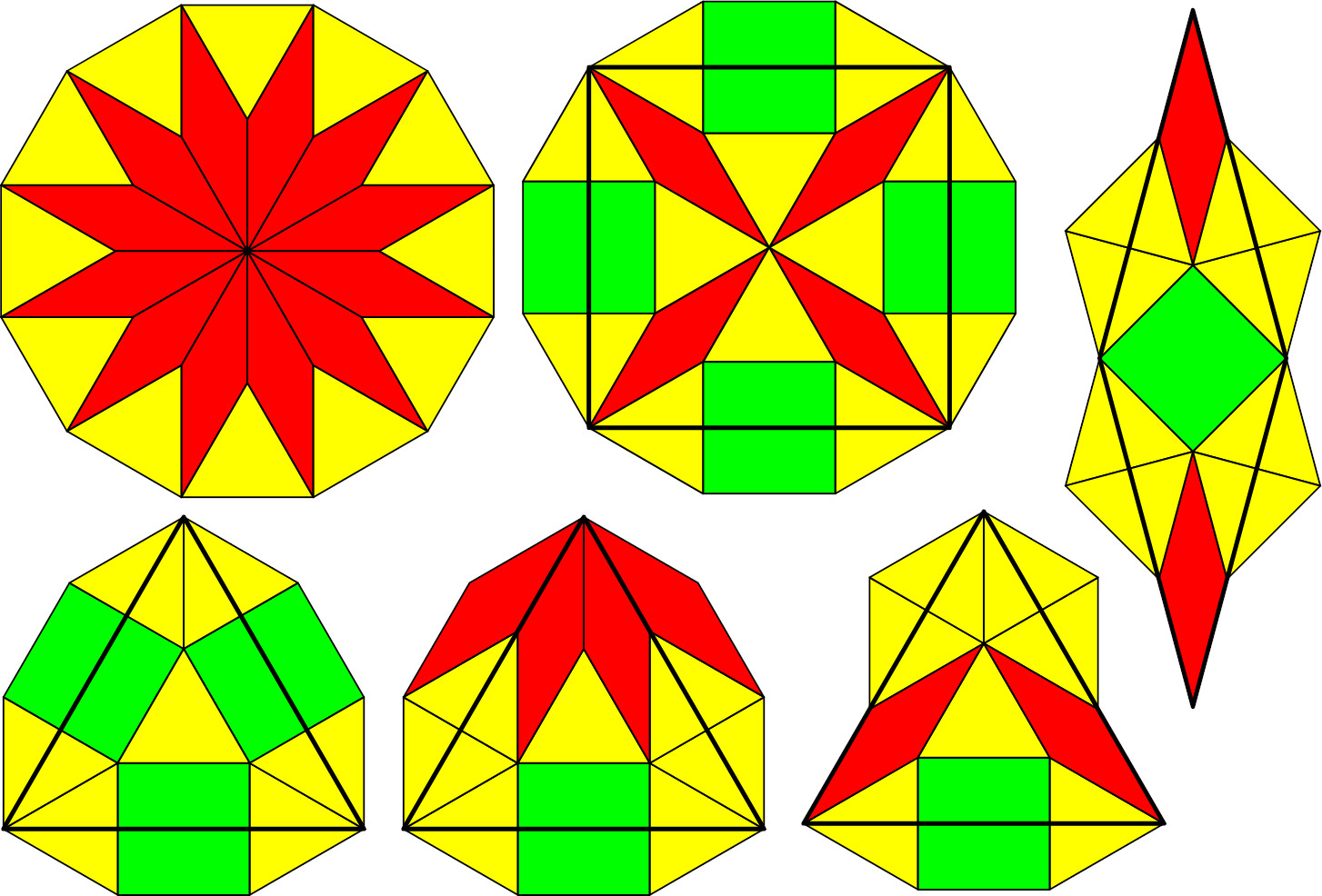}
			\caption{From top left to bottom right: Dissection of a dodecagon into rhombs and triangles. Substitution rules for squares, rhombs and triangles used in the tiling of 1998.}
			\label{fig:theosubstitution}
		}
	\end{minipage}
\end{figure}

\begin{figure}
	\centering
	\begin{minipage}[b]{1\textwidth}
		\textit{
			\centering
			\includegraphics[width=1\linewidth]{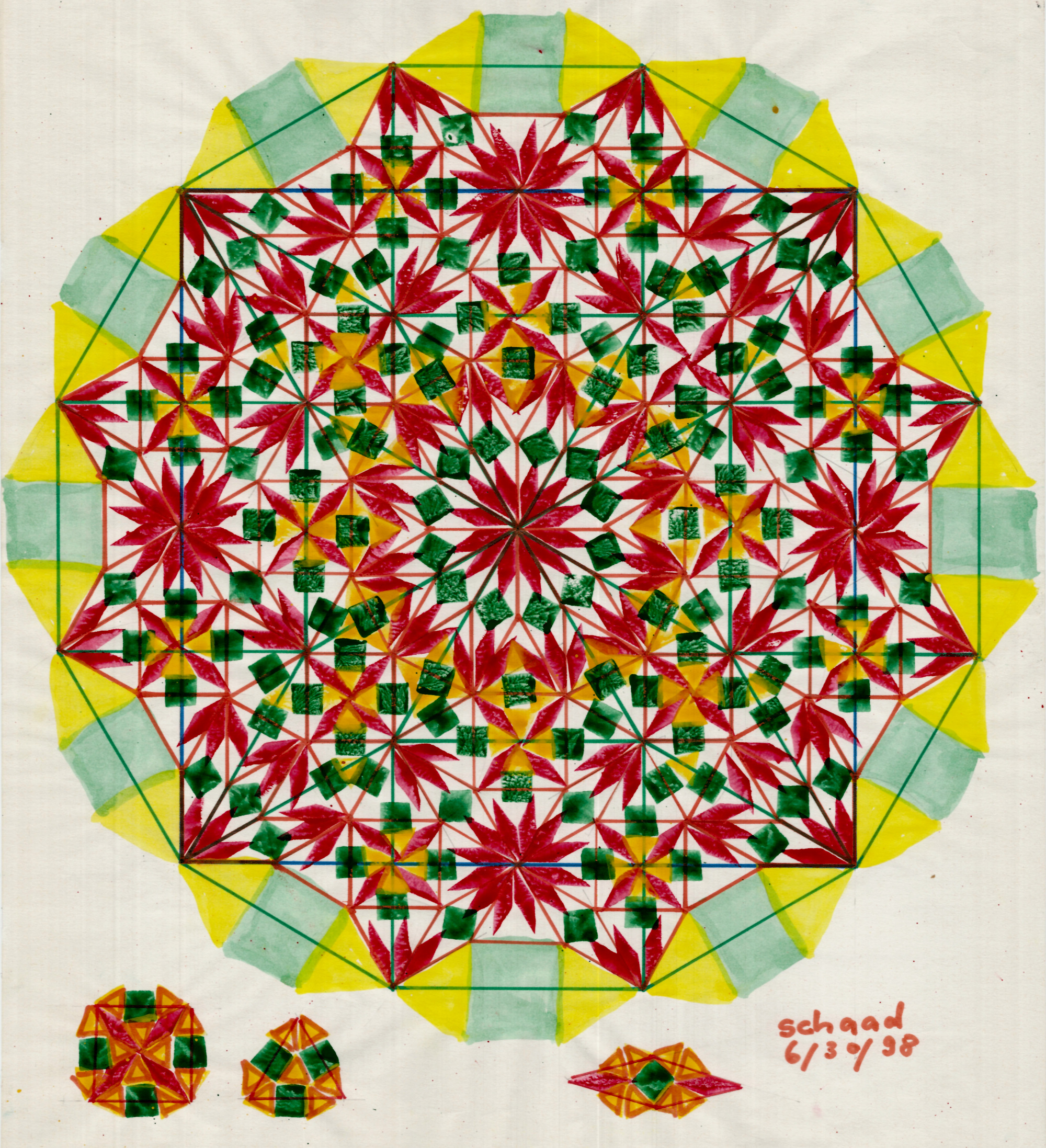}
			\caption{Drawing of 1998: The green outline is the dodecagon of Figure~\ref{fig:theosubstitution}, divided into 4 squares, 4 rhombs, and triangles (green lines). The green squares and rhombs were then further divided (orange lines) in the manner shown at the bottom. The orange shapes were then filled again with squares, rhombs, and triangles using the rules of Figure~\ref{fig:theosubstitution}.}
			\label{fig:tiling1998}
		}
	\end{minipage}
\end{figure}

\section*{The Tiling of 1998}

In 1998, Theo Schaad tried to find substitution rules for a quasiperiodic tiling with 12-fold rotational symmetry. For the ratio between the edges of an inflated tile and the base tiles he chose $1 + \sqrt{3}$ in analogy to the self-similarity ratio $1 + \sqrt{2}$ of the Ammann-Beenker tiling. This is different to the ratio $2 + \sqrt{3}$ of many tilings with 12-fold rotational symmetry and makes it much more difficult to find substitution rules.
He noted that adding 12 equilateral triangles to a rosette of 12 rhombs makes a dodecagon as shown in Figure~\ref{fig:theosubstitution}. Replacing some rhombs with squares and moving triangles as shown in Figure~\ref{fig:replacement}, he obtained a rosette with 4-fold rotational symmetry. It defines a substitution rule for a square with the desired inflation ratio. A substitution rule for rhombs with the same ratio is also easily found as shown in Figure~\ref{fig:theosubstitution}. However, the sides of the inflated square and rhomb are incompatible, and it is impossible to place them adjacent to each other because the triangles do not match. But, luckily, the rhombs and squares do not touch, instead they are connected by equilateral triangles. Three possible substitutions for triangles are shown in Figure~\ref{fig:theosubstitution}, but there exist 10 more variations that are less symmetric. Always on the lookout for symmetric solutions first, the most symmetric one is of 3-fold rotational symmetry with three squares, followed by two mirror-symmetric solutions with two rhombs each and a single square at its base. Actually, only these less symmetric substitutions make it possible to get rosettes of 12-fold rotational symmetry.

To correctly choose triangle substitutions, one must check that the tiles of the next generation match. Theo Schaad accomplished this in Figure~\ref{fig:tiling1998}. First, he drew a large dodecagon and subdivided it with green lines into 4 squares, 4 rhombs, and equilateral triangles. This would be the first generation of substitutions as in Figure~\ref{fig:theosubstitution}. He then repeated the iteration with orange lines. This is now the second generation. Careful examination shows that triangles in an inflated square are of the second type with two rhombs at the top and a square base.
It was now possible to move to the third generation, at least partially, by filling all the orange squares and rhombs with base tiles. Spontaneously, a rosette of orange rhombs had appeared at the center. It was obvious that the triangles that completed the rosette had to be of the third type in Figure~\ref{fig:theosubstitution} with two rhombs pointing towards a square at the base. This filled almost the entire tiling with only a few orange triangles left unfilled. Not surprisingly, Theo Schaad tried to find homes for the 3-fold symmetric triangles with 3 squares. A few could be completed that way, leaving the rest with 2 squares and a rhomb, a half-triangle of type 1 joined by a half-triangle of type 2. Hence, a coherent substitution rule with only the tiles of Figure~\ref{fig:theosubstitution} was not possible. The result was, nevertheless, a large finite tiling of nonperiodic 12-fold rotational symmetry with an inflation ratio of $1 + \sqrt{3}$. It consists of thin rhombs, equilateral triangles and squares. The substitution puzzle was not solved for another 22 years but was greatly helped by the discovery of another tiling by Peter Stampfli in 2012.

\begin{figure}
	\centering
	\begin{minipage}[b]{1\textwidth}
		\textit{
			\centering
			\includegraphics[width=0.9\linewidth]{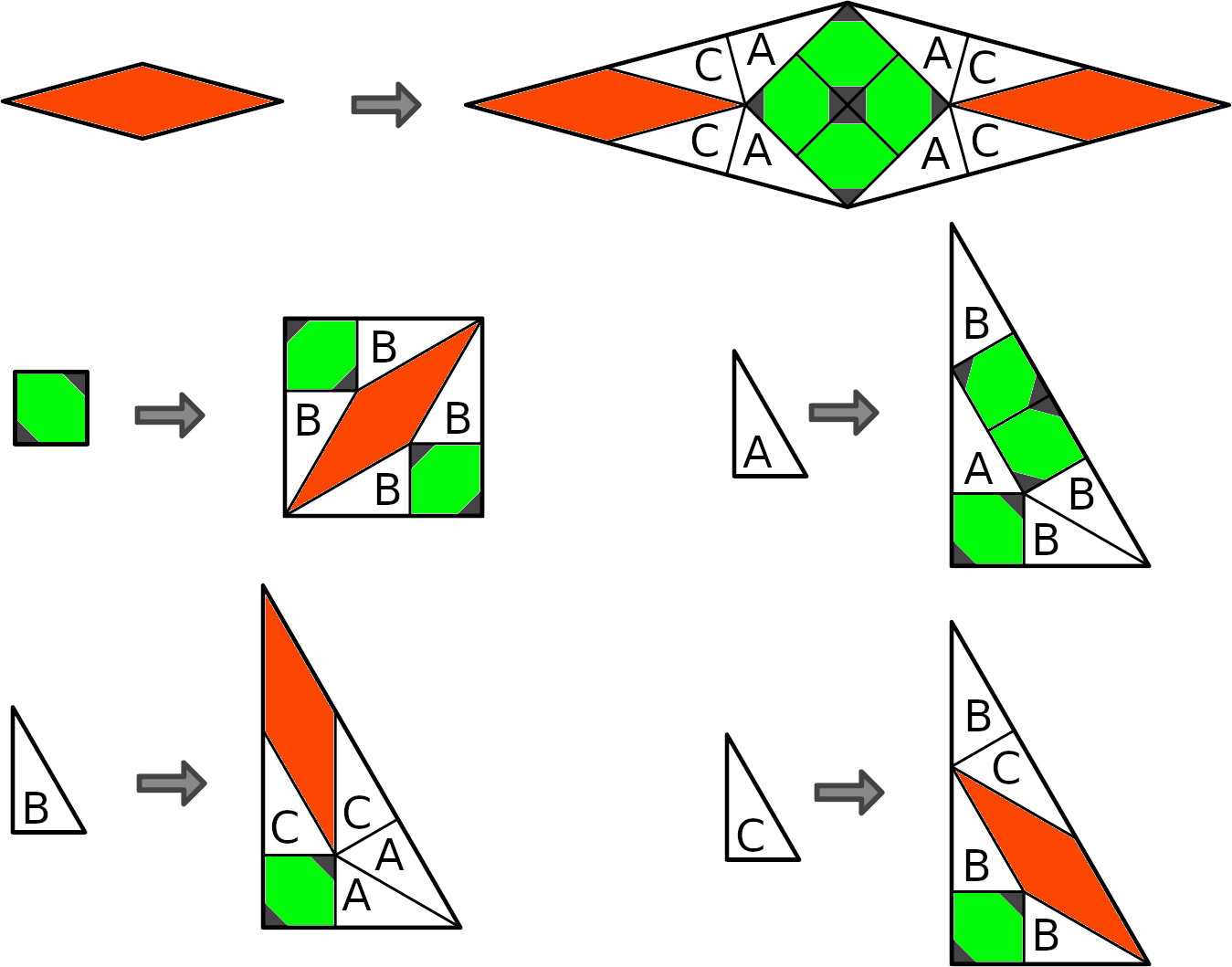}
			\caption{The substitutions of the tiling of 2012. Note the similarity to Figure~\ref{fig:theosubstitution}.}
			\label{fig:peterSubstitution}
		}
	\end{minipage}
\end{figure}

\begin{figure}
	\centering
	\begin{minipage}[b]{1\textwidth}
		\textit{
			\centering
			\includegraphics[width=1\linewidth]{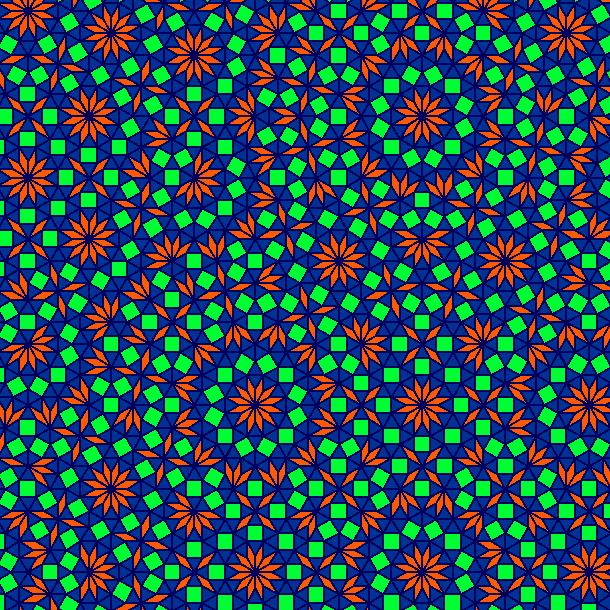}
			\caption{A patch of Peter Stampfli’s 12-fold tiling resulting from the substitutions of Figure~\ref{fig:peterSubstitution}.}
			\label{fig:peter}
		}
	\end{minipage}
\end{figure}

\section*{The Tiling of 2012}

Independently, Peter Stampfli\cite{anotherTiling} created a 12-fold tiling with the same properties in 2012. He encountered similar issues trying to find substitutions for triangles that had to match either an adjacent square or rhomb. The end result was quite similar but not exactly the same. He actually found one of the nine variations presented in the next section.

Since Peter Stampfli worked with a computer program, he could not just fill in triangles as they fit. He had to define an algorithm for a specific iteration scheme. He could then use the program to find out if his scheme made a good tiling. At first, he noted that there are three choices to lay out the edges of an enlarged triangle as shown in Figure~\ref{fig:replacement}. Potentially, there could be 27 different triangles. Some are simply rotations of others, but there are still 13 left. He chose the same substitution rule for squares as in Figure~\ref{fig:theosubstitution} and observed that a triangle next to a square has to have a square and two triangles as substitution at its corresponding side, see the center of Figure~\ref{fig:replacement}. In other cases, when two triangles met, this substitution was also the best symmetric solution. Therefore, every triangle had to have at least one square in its substitution, which left 9 different possibilities. This could be reordered in a much simpler way. He cut the triangles in half and used these half-triangles as base tiles for the substitution scheme. There were only three left and they could be combined in nine different ways to make equilateral triangles. He kept track of them with letters, see Figure~\ref{fig:peterSubstitution}. This simplification had another advantage. By using only a quarter of the previous larger square, he could now inflate all base tiles into exactly similar larger tiles with an inflation ratio of $1 + \sqrt{3}$. The triangles had to be chosen such that half-triangles always met together to give equilateral triangles in the tiling. As Theo Schaad had found too, some substitutions were obvious, e.g., the triangles in the substitution of the squares had to be B triangles (or of the second type in Figure~\ref{fig:theosubstitution}). Likewise, the inflated rhombs had to be filled with C and A triangles. This finally left only very few choices for the triangles. Peter Stampfli found the substitution rules shown in Figure~\ref{fig:peterSubstitution} and created the patch of Figure~\ref{fig:peter}. The result was a potentially infinite quasiperiodic 12-fold tiling with three tile shapes that were themselves composites of smaller base tiles (half-triangles and quarter-squares). Peter Stampfli noted in his blog that the result looked rather complex, joining others that found 12-fold tilings somewhat chaotic. But it was also a choice of presentation. Inflated portions resulting from central 12-fold rosettes or of other rotational symmetries are quite aesthetic.

\begin{figure}
	\centering
	\begin{minipage}[b]{1\textwidth}
		\textit{
			\centering
			\includegraphics[width=0.9\linewidth]{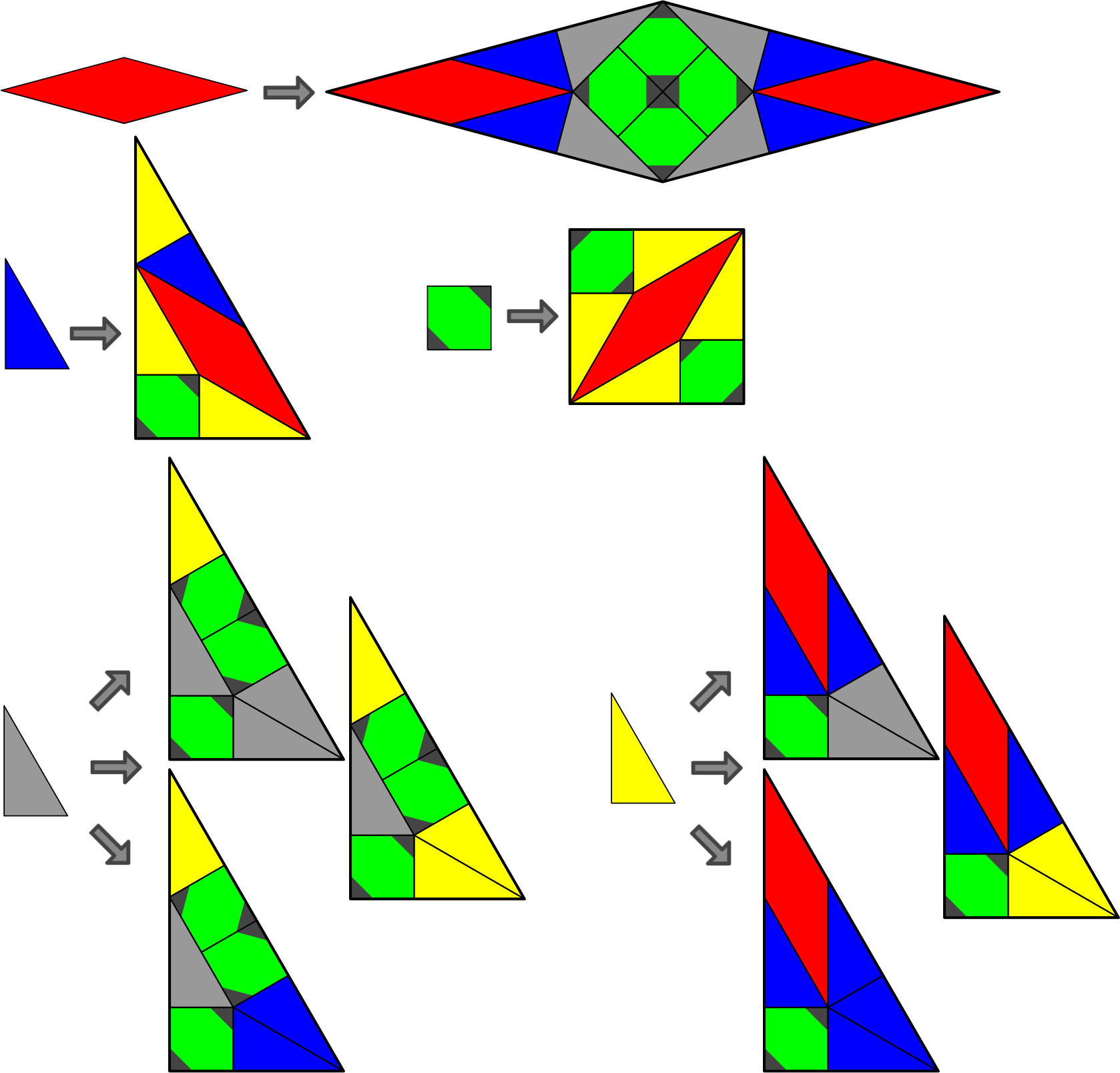}
			\caption{The new substitution rules for a quasiperiodic 12-fold tiling with an inflation factor $1 + \sqrt{3}$. For the gray and yellow triangles we can choose between three substitutions leading to nine different tilings.}
			\label{fig:substitutionFinal}
		}
	\end{minipage}
\end{figure}

\section*{The New Substitution Rules}

\begin{figure}
	\centering
	\begin{minipage}[b]{1\textwidth}
		\textit{
			\centering
			\includegraphics[width=1\linewidth]{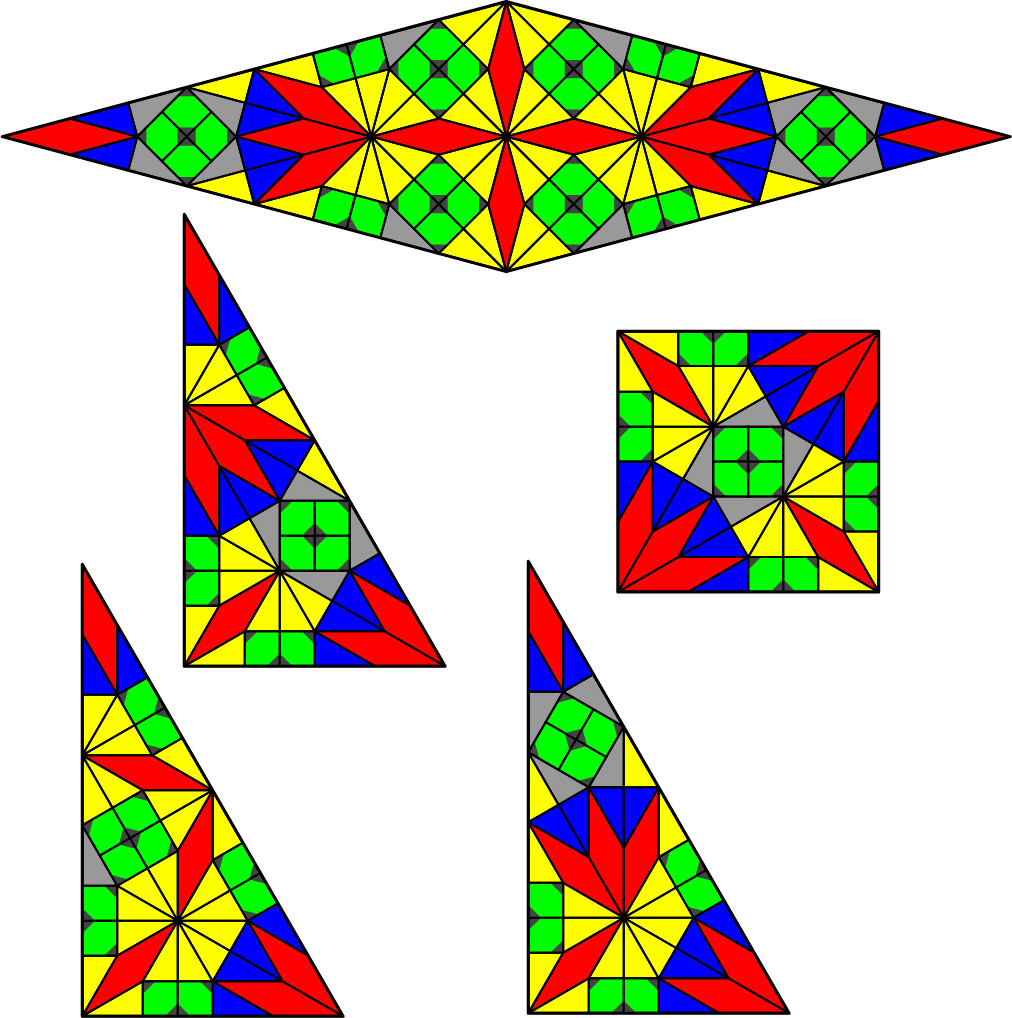}
			\caption{The second generation of one the substitution rules of Figure~\ref{fig:substitutionFinal}. Here, using the second choice in Figure~\ref{fig:substitutionFinal} yellow triangles are put at the $60^\circ$ corner of gray and yellow triangles. This leads to rhombs in the lower corners of the triangles, ultimately generating many rosettes in a larger patch.}
			\label{fig:second}
		}
	\end{minipage}
\end{figure}

We decided to follow Peter Stampfli’s choice of base tiles with a rhomb of unit side length, a square with sides of half unit length, and right triangles with a hypotenuse of unit length. Further research gave the substitution rules of Figure~\ref{fig:substitutionFinal}. Instead of labeling the triangles with letters, we preferred different colors. This makes it considerably easier to follow repeated substitutions as shown in the example given in Figure~\ref{fig:second}. We see that in higher generations two-colored equilateral triangles appear. Here gray and yellow halves are combined. To regain the simplicity of a tiling with only three shapes (rhomb, square and equilateral triangle with unit side lengths) we could draw all triangles using the same color. Comparing the new rules with the tiling of 2012, see Figures~\ref{fig:peterSubstitution} and \ref{fig:peter}, the triangle A corresponds to a gray triangle and follows the second variation, and the triangle B follows the first variation of the yellow triangle. We also discovered that the substitutions of the gray and blue triangles have to use a yellow triangle at the apex because yellow triangles are substituted with rhombs at the apex. This rule showed that there could not be a triangle substitution with only gray triangles as proposed in the 1998 Tiling. One final caveat: whenever two triangles are joined along an edge, it is still possible to exchange at this edge the square and triangles by two rhombs and triangles. However, it would invalidate our goal of making inflated tiles that are similar to the base tiles.

\begin{figure}
	\centering
	\begin{minipage}[b]{1\textwidth}
		\textit{
			\centering
			\includegraphics[width=1\linewidth]{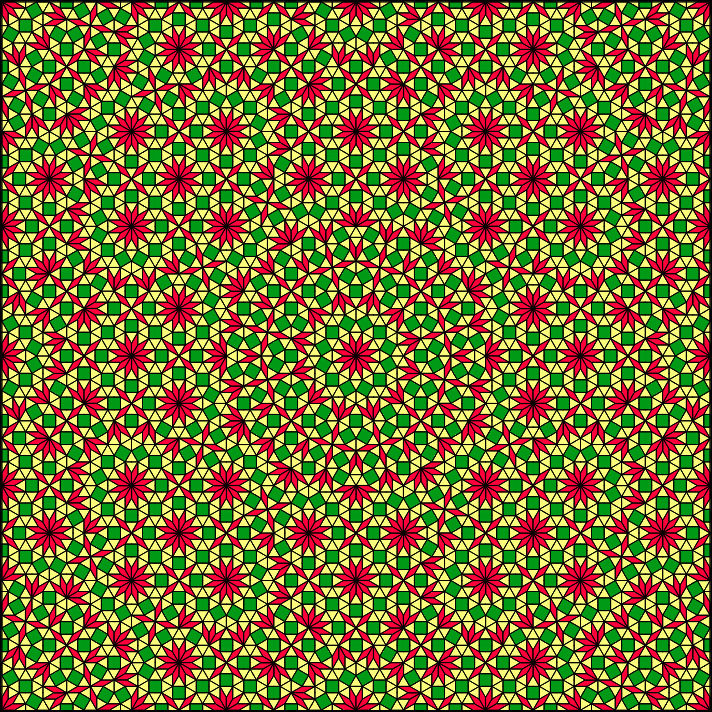}
			\caption{The substitution rules with the choice of Figure~\ref{fig:second} applied four times to a square tile. All triangles are shown in one color. This patch can be compared to the tiling in Figure~\ref{fig:tiling1998}, except that it was iterated to a fourth generation.}
			\label{fig:squareResult}
		}
	\end{minipage}
\end{figure}

\begin{figure}
	\centering
	\begin{minipage}[b]{1\textwidth}
		\textit{
			\centering
			\includegraphics[width=1\linewidth]{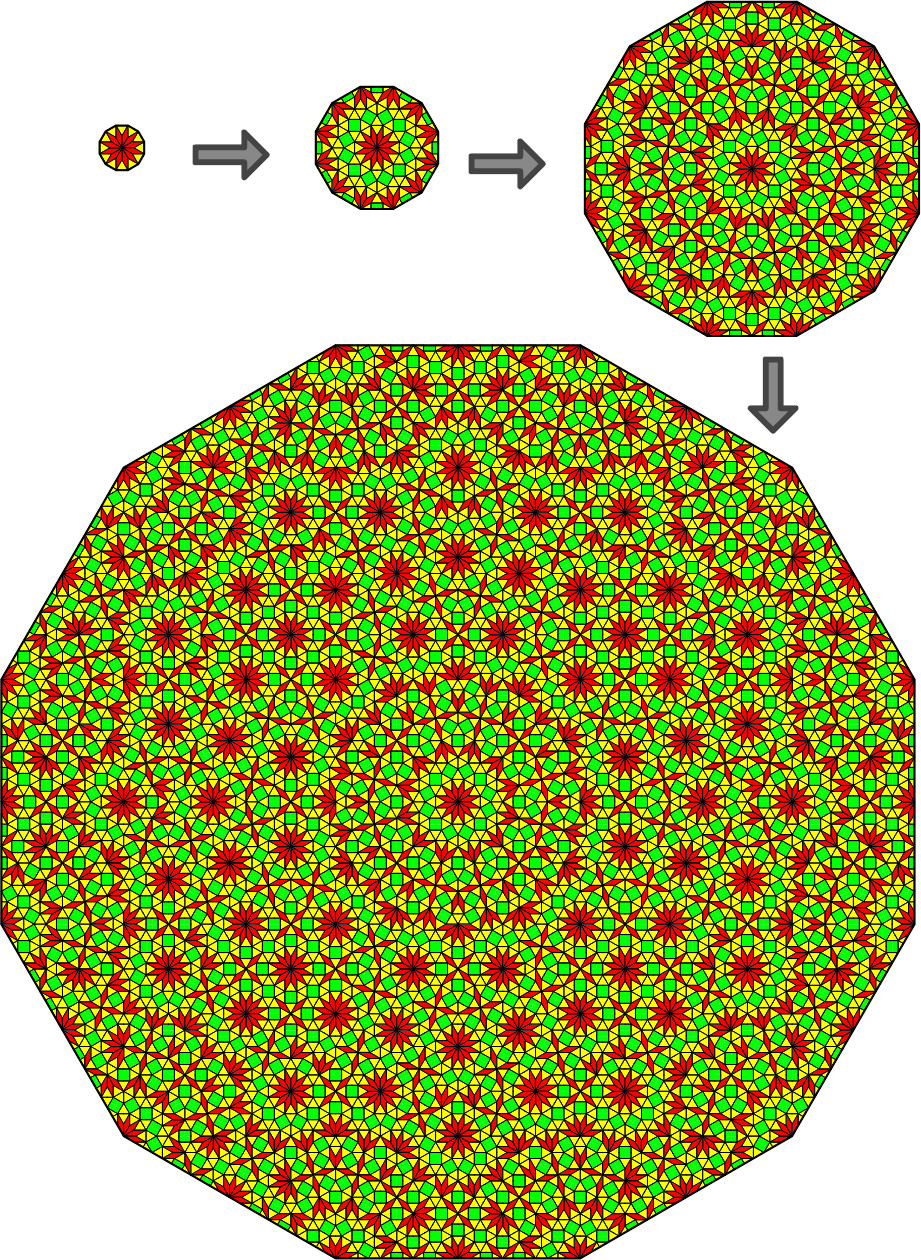}
			\caption{Applying three times the substitution rules with the choice of Figure~\ref{fig:second} to a rosette of twelve rhombs and triangles. This creates a tiling with 12-fold rotational symmetry.}
			\label{fig:dodecas}
		}
	\end{minipage}
\end{figure}

\section*{Results}

The reader can explore this tiling in his browser using a computer program written by P.~Stampfli\cite{app}. Here, we present some of its results. If we use only one color for the three triangle substitutions, the resulting equilateral triangles will be indistinguishable from each other and the final tiling will have only three apparent base shapes: square, rhomb, and equilateral triangle. Figure~\ref{fig:squareResult} shows such an example by applying the substitution rules to a square four times. This replicates the steps Theo Schaad took in his original tiling except that it goes one generation further. In Figure~\ref{fig:dodecas}, we apply the substitution rules three times to the rosette originally shown in Figure~\ref{fig:theosubstitution}. It consists of 12 rhombs and equilateral triangles. Each substitution of a rhomb contains smaller rhombs at its acute corner, thus the rhomb itself becomes part of a tiling of larger rhombs. This process can be thought of as adding new tiles to the original rhomb.  In the case of the enlarged rosette, the rosette is again mapped into itself leaving the original rosette at its center. We can look at the substitution as adding new tiles around the original rosette. Each addition is reminiscent of the nesting sequence of Russian Matryoshka puppets.

\section*{Summary}

We presented a quasiperiodic 12-fold tiling of squares, equal sided triangles, and 30 degree rhombs. It is defined and built by its substitution rules. An inflated patch of the tiling does not change its shape. Only its size is increased by the self-similarity ratio of $1 + \sqrt{3}$. The tiling is edge-to-edge without gaps or overlaps and the substitution process can be repeated to make an arbitrarily large tiling. 
The authors would like to thank Robert Ingalls, known for his work on octagonal and decagonal tilings, for providing the source material of 12-fold patterns by J. Socolar and D. Haussler et al.
% -----------------------------------------------------------------------------
% references

{\setlength{\baselineskip}{13pt} % tighten line spacing for bibliography
	\raggedright				% no right justification for References
	
}

% -----------------------------------------------------------------------------


\begin{thebibliography}{99}
		
		\bibitem{app}
		\url{http://geometricolor.ch/twelveOneSqrt3.html}
		
		\bibitem{encyclopedia}
		D. Frettlöh, E. Harriss and F. Gähler, \textit{Tilings encyclopedia}, \url{https://tilings.math.uni-bielefeld.de/}.
		
		\bibitem{stampfli2}
		\url{https://geometricolor.wordpress.com/2012/09/29/finding-an-iteration-method-for-the-stampfli-tiling-mission-impossible/}
		
		\bibitem{
			ammannBeenker}
		\url{https://tilings.math.uni-bielefeld.de/substitution/ammann-beenker/}
		
		\bibitem{rorschach}
		\url{https://tilings.math.uni-bielefeld.de/substitution/rorschach/}
		
		\bibitem{hofstetter}
		\url{https://tilings.math.uni-bielefeld.de/substitution/hofstetter-4fd-arrows/}
		
		\bibitem{anotherTiling}
		P. Stampfli,  \textit{Another tiling of dodecagonal symmetry},
		\url{https://geometricolor.wordpress.com/2012/07/29/another-tiling-of-dodecagonal-symmetry}
		
		\bibitem{stuererDeludy}
		W. Steurer and S. Deloudi, \textit{Crystallography	of Quasicrystals, Concepts, Methods and Structures}, 2009, Springer, Heidelberg
	
		\bibitem{geometricolor}
		P. Stampfli, \textit{Geometry in color}, \url{https://geometricolor.wordpress.com/}
		
		\bibitem{gruenbaum}
		B. Grünbaum and G. Shephard, \textit{Tilings and Patterns}, W.H. Freeman, 1987

		\bibitem{socolarOctagonal} J. Socolar, \textit{Simple octagonal and dodecagonal quasicrystals}, Phys. Rev. B 39, (1989), 10519	
		
		\bibitem{periodicPointSet}
		M. Baake, R. Moody and  M. Schlottmann, \textit{Limit-(quasi) periodic point sets as quasicrystals with p-adic internal spaces}, J. Phys. A: Math. Gen. 31, (1998), 5755-5765
		
		\bibitem{watanabe}
		Y.Watanabe, T. Soma, and M. Ito, \textit{A new quasiperiodic tiling with dodecagonal symmetry},
		Acta Crystallographica Section A: Foundations and Advances, Volume 51 (6) – Jan 1, 1995
		
		\bibitem{gaehlerCrstallography}
		F. Gähler, \textit{Crystallography of dodecagonal quasicrystals} in \textit{Quasicrystalline materials : Proceedings of the I.L.L. / Codest Workshop, Grenoble, 21 - 25 March 1988}, World Scientific, Singapore, (1988), 272-284 
		
		\bibitem{singThesis}
		B. Sing, \textit{Beyond Pisot Substitutions}, PhD thesis, Univ. Bielefeld, 2007
		
		\bibitem{sololarThesis}
		J. Socolar, \textit{Quasilattices and Quasicrystals}, PhD Thesis, University of Pennsylvania, 1987
		
		\bibitem{stampfliHelvPhys}
		P. Stampfli, Helv. Phys. Acta 59, 1260 (1986)
		
		\bibitem{altourian}
		V. Altourian and P. Stampfli, \textit{Cover Stories: Making the Graphene Quasicrystal Cover}, Science, Aug 2018

		\bibitem{hausslerNissen}
		D. Haussler, H.U. Nissen, and R. Lueck, \textit{Dodecagonal Tilings Derived as Duals from Quasiperiodic Ammann-Grids}, phys. stat. sol. (a) 146, (1994), 425
				
	\end{thebibliography}
\end{document}